\documentclass[12pt]{article}
\usepackage[english]{babel}
\usepackage{amsfonts}
\usepackage[dvips]{graphicx}

\begin{document}

\title{\bf Searching optimal shapes for blades of a fan}

\date{March 2008}

\author{\bf Gianluca Argentini \\
\normalsize gianluca.argentini@riellogroup.com \\
\normalsize gianluca.argentini@gmail.com \\
\textit{Research \& Development Department}\\
\textit{Riello Burners}, 37045 San Pietro di Legnago (Verona), Italy}

\maketitle

\begin{abstract}
A nonlinear differential equation about optimal shapes for blades of a fan. A boundary value differential problem from engineering, geometrical or physical bonds. A relation between linear profiles and constant speed along the side under flow.\\

\end{abstract}

\section{The differential equation}

Consider the 2D axial projection of a centrifugal forward fan (\cite{bleier}). Air is sucked from a central zone and channeled between two blades. Then, centrifugal force put it towards the external zone of the fan. From the point of view of a single blade, air flow reaches its geometrical basis and go through the blades channel. In this work the fundamental hypothesis is that {\it the flow field is tangent to the profile of a blade at its initial point}.

\begin{figure}[ht]\label{plotFig1}
	\begin{center}
	\includegraphics[width=8cm]{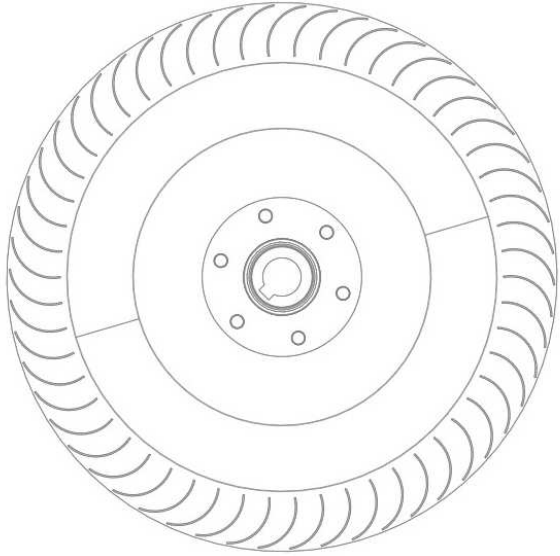}
	\caption{\it{A forward fan's impeller: the concave side of the blades receives the fluid flow.}}
	\end{center}
\end{figure}

It can be shown that this condition is suitable in the case of noise reduction. But what about the profile of the blade, in the case of geometrical or engineering bonds?\\

\begin{figure}[ht]\label{plotFig2}
	\begin{center}
	\includegraphics[width=6cm]{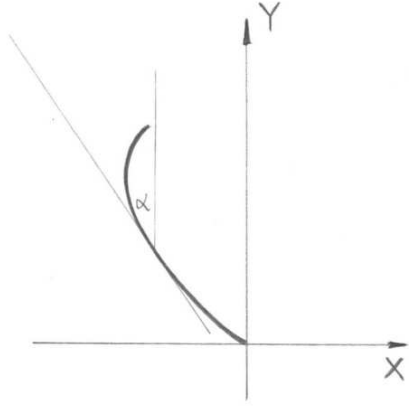}
	\caption{\it{The reference cartesian system local to the blade.}}
	\end{center}
\end{figure}

Let \{$X,Y$\} a cartesian system local to a blade, with $Y$-axis containing the blade's initial geometrical point (see Fig.2). Assume that the blade's profile is the graphic of a smooth function $X=F(Y)$. Also we assume that {\it the flow velocity field is tangent along the concave side}, that is we neglect physical effects as boundary layer or flow separation (\cite{rouse}). Then, the acceleration along the blade's side is the sum of the tangential and the centripetal accelerations, that is (see \cite{courant})

\begin{equation}\label{accelerations}
	\mathbf{a}_t = \dot{v}\hspace{0.1cm}(sin\alpha,cos\alpha), \hspace{1cm} \mathbf{a}_c = \frac{v^2}{r_c}(cos\alpha,sin\alpha)
\end{equation}

\noindent where $v$ is the local speed in the \{$X,Y$\} system, $v^2 = \dot{X}^2+\dot{Y}^2$, and $r_c$ is the radius of curvature:

\begin{equation}\label{radiusCurvature}
	r_c = \frac{(1+F_Y^2)^{\frac{3}{2}}}{|F_{YY}|}
\end{equation}

\noindent Note that we consider $\alpha \geq 0$ if $F_Y \leq 0$, negative otherwise. Then the $X$-component of the total local acceleration is

\begin{equation}\label{XlocalAcc}
	\ddot{X} = \dot{v}\hspace{0.1cm}sin\hspace{0.1cm}\alpha + \frac{v^2}{r_c}cos\hspace{0.1cm}\alpha
\end{equation}

\noindent Now we try to write previous equation using only $Y$ and its derivatives. First of all, as $X(t) = F(Y(t))$, by chain rule we have

\begin{equation}\label{dotX}
	\dot{X} = F_Y\hspace{0.1cm}\dot{Y}
\end{equation}

\noindent

\begin{equation}
	\ddot{X} = F_{YY}\dot{Y}^2 + F_Y\ddot{Y}
\end{equation}

\noindent so that

\begin{equation}\label{speed}
	v^2 = \dot{X}^2 + \dot{Y}^2 = (1+F_Y^2)\dot{Y}^2
\end{equation}

\noindent and

\begin{equation}\label{dotspeed}
	2\hspace{0.1cm}v\hspace{0.1cm}\dot{v} = 2F_YF_{YY}\dot{Y}^3 + 2(1+F_Y^2)\dot{Y}\ddot{Y}
\end{equation}

\noindent Also, note that $tg\hspace{0.1cm}\alpha = - F_Y$, so using usual trigonometric formulas we can write

\begin{eqnarray}\label{sincos}
	sin\alpha = \frac{tg\alpha}{\sqrt{1+tg^2\alpha}} = - \frac{F_Y}{\sqrt{1+F_Y^2}}\\
	cos\alpha = \frac{1}{\sqrt{1+tg^2\alpha}} = \frac{1}{\sqrt{1+F_Y^2}} \nonumber
\end{eqnarray}
	
\noindent With all these previous identities and the assumption $F_{YY} > 0$, after some simplifications equation (\ref{XlocalAcc}) becomes

\begin{equation}\label{diffEq1}
	F_Y(1+F_Y^2)\ddot{Y} +F_Y^2F_{YY}\dot{Y}^2 = 0
\end{equation}

\noindent In the points of the blade where $F_Y = 0$ this equation is banally satisfied, so we can say that in the points where $F_Y$ is not null the blade equation is

\begin{equation}\label{diffEq2}
	(1+F_Y^2)\ddot{Y} +F_YF_{YY}\dot{Y}^2 = 0
\end{equation}

\noindent Previous is a nonlinear differential equation that can be used in two manners: if one assigns the blade profile, that is the function $F(Y)$, then by (\ref{diffEq2}) it can be computed the law motion $Y=Y(t)$ along the blade side; if one wants to assigne a particular motion $Y(t)$, then it can compute the blade profiles $F(Y)$ which satisfy that law.

\section{A differential problem}

Suppose that engineering project or physical bonds require a particular law motion $Y=Y(t)$ such that $\dot{Y}=g(Y)$. Then $\ddot{Y}=g_Y(Y)g(Y)$. Note that, for our geometrical assumptions, we have a boundary condition $F(b)=0$, where $(0,b)$ are the cartesian coordinates of the blade's initial point. Also, from physical assumptions about the tangential flow at blade's initial point, we have a bondary condition on the derivative: $F_Y(b)=m_0$, being $m_0$ the slope of the flow field at impeller inlet. Then the blade profile is determined by the differential problem

\begin{equation}\label{diffProblem1}
\left\{
\begin{array}{lll}
	\nonumber g_Y(Y)g(Y)(1+F_Y^2) + g^2(Y)F_YF_{YY} = 0\\
	F(b) = 0\\
	\nonumber F_Y(b) = m_0
\end{array}
\right.
\end{equation}

\noindent A first simple but useful consideration is the following. Suppose to have blades with a null curvature, that is $F_{YY}=0$ for every $Y$. Then, from (\ref{diffProblem1}) follows 

\begin{equation}
	g_Y(Y)g(Y)(1+F_Y^2)=0
\end{equation}

\noindent Therefore $g_Y=0$ for every $Y$, that is $Y(t)=at+b$ being $\ddot{Y}=0$. Vice versa, if $Y(t)$ is linear on $t$ variable, then $g_Y=0$ for every $Y$, and the differential equation becomes 

\begin{equation}
  g^2(Y)F_YF_{YY} = 0
\end{equation}

\noindent which gives $F_{YY}=0$. Note that, in this case, from (\ref{dotX}) follows that $X(t)$ is linear too. Therefore, {\it the only profiles with constant speed along the downwind side are linear blades}.

\end{document}